\documentclass[oneside,12pt]{amsart}
\usepackage{amsmath,amsfonts, latexsym,amssymb}
\usepackage[active]{srcltx}

\theoremstyle{definition}

\numberwithin{equation}{section}


\begin{document}

\title[]{A short survey  of the game Bulgarian solitaire and related games}

\author{Romeo Me\v strovi\' c}

\address{Maritime Faculty Kotor, University of Montenegro, Dobrota,
 85330 Kotor, Montenegro} \email{romeo@ucg.ac.me}

  \begin{abstract}
Let $N$ be an arbitrary  positive integer and let 
$\lambda=(\lambda_1, \lambda_2,\ldots,\lambda_l)$ be a partition of $N$
of length $l$,  i.e., $\sum_{i=1}^{l}\lambda_i=N$
with parts $\lambda_1\ge \lambda_2\ge\ldots\ge \lambda_l\ge 1$.
 Define $T(\lambda)$ as the partition of $N$ with parts 
$l,\lambda_1-1, \lambda_2-1,\ldots,\lambda_l-1$, ignoring any
 zeros that might occur.

Starting with a partition $\lambda$ of $N$,
we describe Bulgarian solitaire by repeatedly applying the shift
operation $T$ to obtain the sequence of partitions
     $$
     \lambda, T(\lambda),T^2(\lambda),\ldots .
      $$
  We say a partition $\mu$ of $N$ is $T$-{\it cyclic} if $T^i(\mu)=\mu$
  for some $i\ge 1$. In 1982 Brandt \cite{br} characterized all  $T$-{\it cyclic}
  partitions for Bulgarian solitaire. Bulgarian solitaire is a dynamical 
system on integer partition of a positive integer $n$ which converges 
to a unique fixed point if $N=1+2+\cdots +k$ is a triangular number.

In this paper we give a  short survey  of the game Bulgarian solitaire
 and several variations of this  game.
   \end{abstract}

\maketitle

\noindent 2010 {\it Mathematics Subject Classification}: 
 Primary 05A17,  01A60, Secondary 11P81, 91A46.

\vspace*{1mm}

\noindent{\it Keywords and phrases}: 
Bulgarian solitaire, cycle of partitions,
shift operation, triangular number, partition of the integer. 
 \section{The game Bulgarian solitaire}

\vspace{2mm}

The following popular mathematical card game,
popularized (but not introduced) by {\it blast} Martin Gardner in 1983 
\cite{ga}, is called {\it Bulgarian solitaire}.\\

{\bf The game Bulgarian solitaire}. {\it Initially, we are given $N$ cards disposed in several piles. 
A move consists of removing exactly one card from each pile and forming a 
new pile. The operation is repeated over and over.}\\

If the number of cards $N$ is a triangular number, i.e., $N$ 
is of the form $N=1+2+\cdots +k=k(k+1)/2$
for some positive integer  $k$, a remarkable fact is that, starting from any initial
configuration, after a finite number of moves the Bulgarian solitaire
will reach the stable configuration formed by piles of sizes $1,2,\ldots, k$.
This result was proved in 1982 by  Brandt (\cite{br}, the assertion after the proof
of Theorem 4, p. 484). If $N$ is a triangular number, then this game 
is also called {\it Karatsuba solitaire} (\cite{hkk} and \cite{man}).

It was  also considered in \cite{br} the case when the
number of cards
is not triangular. Since a deck has only finitely many layouts, the game
of  Bulgarian solitaire must cycle.                                 
Brandt characterized and counted all cycles for any given
deck  size \cite[Theorem 5]{br}.  

A survey of the earlier history 
of the game Bulgarian solitaire, including its name and summarizing subsequent 
research, was  presented in 2012 by B. Hopkins \cite{ho}. 
Hopkins pointed out that this game was popularized by M. Gardner in 
1983 \cite{ga} with the unusual name {\it Bulgarian solitaire}.

However,  in \cite{ho} Hopkins reported  that this game was  
discovered in 1980 by Konstantin Oskolkov by demonstrating the 
particular case of the game Bulgarian solitaire  concerning five piles with heights
$3,1,4,1,6$, namely, he wrote: 

{\it ``Around $1980$, Konstantin Oskolkov of the Steklov Mathematical Institute
in Moscow traveled by train to give a talk in Leningrad 
{\rm (now Saint Petesburg} {\rm )}. A man on the train told him
of the problem, although the details of the dialog above are fictional. 
Oskolkov shared this with his colleagues at the institute; reportedly
when one number theorist heard about it ``his face a Satanic expression",
he ran to his office, closed the door and did not come out until he 
solved the problem''}.  
 
Notice that problem about {\it BS} firstly appeared in 1981 in the famous 
Russian journal 
{\it Kvant} \cite{to}. Andrei Toom, a research scientist at Moscow State 
University, published a related solution in 1981 \cite{to}.
Notice that his proof was also included in  1981  Book on Mathematical 
Olympiads \cite{vgrt}. In the same year, the solution of $BS$
was published  by B. Bojanov  in a Bulgarian high school mathematics 
journal  \cite{bo} (here {\it BS} was defined in the context of heaps of 
balls) (see \cite{dr}).  The solution   of {\it BS} was also published in 1981 by 
Eriksson in Sweden mathematical journal \cite{er}. Hopkins \cite{ho}
noticed that the three  1981 solutions of {\it BS} use 
very similar ideas to show that, for $n=T_k$    there is indeeed 
just the single fixed point $\tau_k=(k,k-1,\ldots,2,1)$.

Another two proofs about Bulgarian solitaire were established in 1982 by  Brandt
\cite{br} (a proof in the setting of integer partitions),
and in 1985 by Akin and  Davis \cite{ad}. 
 An induction proof of a result about Bulgarian solitaire was given in 2010
by the author of this paper \cite{me1}.    
Using P\'{o}lya enumeration theory (see \cite{bi}), Brandt \cite{br}
 derived the related formula.
Another combinatorial proof of this formula was  given by the author of 
of this paper in \cite{me3}.

D. Knuth \cite{hk} started a fall 1982 programming and problem-solving seminar 
with {\it BS}. In 1983 {\it blast} M. Gardner \cite{ga} popularized 
{\it BS} in his mathematics column from {\it Scientific American}.
D. Knuth has also conjectured that for triangular number $N$ the length 
of the game (i.e., the number of moves before the final position is reached)
is at most $k(k-1)$. This conjecture was solved 
affirmatevly 
in 1985 by K. Igusa  \cite{ig}. The proof of a generalization 
of this conjecture for arbitrary $N$ was established in 1991 by G. Etienne
\cite[Theorem 5.1]{et} (the paper submitted in 1984!). Etienne also 
proposed related additional conjectures \cite[Conjecture 5.2]{et}
 
 In 1987 H.-J. Bentz gave another solution of Bulgarian solitaire \cite{ben}.

Notice that from the proof of a result about Bulgarian solitaire 
 given in  \cite{me1} 
 may be easily deduced the solution of the following interesting 
counting problem on integer sequences.

\vspace{2mm}

{\bf Problem.} {\it Denote by $\# S$ the cardinality of a finite sequence $S$. 
Let $k$ be any fixed nonnegative number, and let  
$(a_n)_{n=1}^{\infty}$ be a sequence of integers $($not necessarily 
nonnegative$)$
satisfying the following conditions:

$\,\,(i)$ the first $k$ terms of a sequence $(a_n)_{n=1}^{\infty}$ are arbitrary, and

$(ii)$ for each $n>k$ there holds
  $$
a_n=\# \{i:\,1\le i\le n-1\,\, {\rm and}\,\, a_i+i<n \}.
  $$ 
Show that $(a_n)_{n=1}^{\infty}$ is a periodic sequence.}\\

The solution of the above problem was presented in 1981 by B. Bojanov 
\cite{bo}.

In  \cite{kkh} Mizan and Khan (also see Abstract of this paper
available at \cite{kh}), the authors discussed/revisted  Andrei Toom's proof 
of Bulgarian solitaire that appeared in 1981 in {\it Kvant}, 
and show how an application of the Chinese Remainder Theorem
allows us to generalize the proof.

Notice that our list of references includes all twenty items 
from \cite{ho}, as well as 
14 other references and sources.

\section{A formal approach to the mathemacical card game Bulgarian solitaire}

The foundations of the theory of integer partitions were laid by Leonhard 
Euler. A goood introduction to this subject are the books \cite{a} and 
 \cite{ae}. 
Let us now define the game  formally.
  Let $N$ be a positive integer and let $\lambda$ be a partition of
 $N$ having $l$ parts written
 $(\lambda_1, \lambda_2,\ldots,\lambda_l)$ in non-increasing order;
that is, $N=\lambda_1+ \lambda_2+\cdots+\lambda_l$ with positive
integers $\lambda_1\ge \lambda_2\ge\ldots\ge \lambda_l\ge 1$.
 Define $T(\lambda)$ as the partition of $n$ with parts $\lambda_1-1, \lambda_2-1,\ldots,\lambda_l-1, l$, ignoring any
 zeros that might occur. So $T^i(\lambda)$ ($i=1,2,\ldots$) denotes the
 partition obtained by successively applying the shift operation $T$ to
 $\lambda$ a total of $i$ times.
 
  Starting with a partition $\lambda$,
 we describe Bulgarian solitaire by repeatedly applying the shift
 operation to obtain the sequence of partitions
     $$
     \lambda, T(\lambda),T^2(\lambda),\ldots .
      $$
    We say a partition $\mu$ of $N$ is $T$-{\it cyclic} if $T^i(\mu)=\mu$
  for some $i\ge 1$. 

    If $N$ is arbitrary, Brandt noted  that
 repeated application of $T$ leads into a cycle of partitions, since 
 there are only  a finite number of these. 
Furthermore,  a cycle
 of  partitions is completely determined by the sequence
 of the consecutive lengths of the partitions in the cycle.
 Motivated by this fact, Brandt (\cite[p. 483]{br}[2] defines the set 
$M_n$ by\\
 
\noindent $
(1)\qquad\qquad\qquad\qquad  M_n=\{\sigma =(\sigma_i)_{i\in{\mathbf Z}}:\,\max\sigma_i=n,\,\,$ 

\noindent\qquad\qquad\qquad\qquad\quad ${\rm where
\, \, for\,\, all \,\,} i,\,
\,
 \sigma_i=|\{\sigma_j|j<i,\sigma_j\ge i-j\,\,
\}|\},
   $\\
   
\noindent where $|S|$ denotes the cardinality of a set $S$.
    If $\sigma\in M_n$, then by Proposition 2 in \cite{br},
   $\sigma_i\in \{n, n-1\}$ for all $i\in{\mathbf Z}$.
   As an easy consequence of this
   fact, Brandt (cf. proof of Theorem 5 in \cite{br}; 
]also see       \cite{ad}, Theorems 4 and 5,  \cite{gh}, Theorem 2.1 and
      Etienne   \cite{et}),
   characterized all $T$-cyclic partitions for an arbitrary positive 
integer $N$.
   This result is given as follows. 
 
   \vspace{2mm}

 \noindent{\bf Theorem.\cite{et}}
   {\it Let  $N=1+2+\cdots +k +r$,
   $0\le r\le k$. Then a partition $\lambda$ of $N$ is  $T$-cyclic if and
   only if $\lambda$ has the form
     $$
     (k+\delta_k, k-1+\delta_{k-1}, \ldots, 1+\delta_1,\delta_0),
      $$
      where each $\delta_i$ is $0$ or $1$ and $\sum_{i=0}^{k}\delta_i=r$.}

 \vspace{2mm}

In particular (see the assertion after the proof  of Theorem 4 in \cite{br}).
For a triangular number $N$ we obtain the following result
quoted by Gardner in 1983 \cite{ga}--{\it Brandt's Equilibrium Theorem}. 

   \vspace{2mm}
 \noindent{\bf Corollary} (Brandt's Equilibrium Theorem)  
{\it If  $N=1+2+\cdots +k$,
then $(k,k-1,\ldots,1)$ is  the unique  $T$-cyclic partition of $N$.}

 \vspace{2mm}

 Recall that the above theorem follows from Theorem 4 in 
\cite{ad}  whose proof is based on Brandt's result. Theorem 5 in \cite{ad}
 which  is proved  directly, also gives a  description of all $T$-cyclic 

partitions for
 arbitrary $N$ as in above theorem.  The above corollary is proved
 by Etienne \cite{et} by introducing a natural array representation
 of a partition $\lambda$. The idea in his proof is applied
 in the proof of Theorem 2.1 in \cite{gh} 
(the above theorem) to general  $N$.

\section{Some variations and new variants of the card game Bulgarian solitaire}
For information about the earlier history of the game Bulgarian solitaire 
and a summary of subsequent research, see reviews by Hopkins \cite{ho}
and Drensky \cite{dr}.

Many variants of Bulgarian solitaire have been suggested in the literarure 
(see \cite{dr} for an extensive survey,) such as in 1985 
{\it Austrian solitaire} by E. Akin and M. Davis
 \cite{ad} (also see \cite{bas}), in  1992 {\it Montreal solitaire} 
by C. Cannings and J. Haigh (England) \cite{ch}, in 1995  
by R. Servedio and Y.N. Yeh \cite{sy}, in 1995 (without name) by Y.N. Yeh,
in 1997 {\it Carolina solitaire} by Andrej Andreev 
(Bulgaria) \cite{man}), \cite{ch}, in 2003 {\it Random Bulgarian solitaire}
by  Popov  \cite{po}, in 2004 
{\it Two-handred Bulgarian solitaire} by Tim  Bancroff  \cite{man}
and in 2012 (the game without name) by B. Hopkins  \cite{jo}.
The PhD thesis \cite{jo} deals with processes on integer partitions and their 
limit shapes, with focus on deterministic and stochactic variants on Bulgarian 
solitaire.  In 2016 Olson \cite{ol} presented a generalization 
of Bulgarian solitaire, the so-called $\sigma$-{\it Bulgarian solitaire},
in which multiple cards may be picked from a single pile.
For another generalization of Bulgarian solitaire, see \cite{ejs}.

Bulgarian solitaire and its variants are extensively researched in 
{\it Combinatorial Game Theory}. S. Dor\'{e}e (of Augsburg College, 
Minneapolis, USA) called {\it BS} ``{\it a somewhat distant 
relative of two-player African pebble games Mancala}'' (see \cite{man}
and \cite{d}).

In \cite{hn} Harris and Nguyen introduced a new representation of Bulgarian 
solitaire that is convenient for the study of related generating functions
(also see \cite{ej}).
In \cite{hn} it was also  proved two instances 
of Pham's conjecture \cite{ph}.

\section{Conclusion}
In this   article,  we present a survey of variations 
of  Bulgarian solitaire game. All these games are inspired by the 
popular mathematical card game Bulgarian solitaire. 
The popularity of the Bulgarian solitaire started  around 1980. 
Note that  in 1983 the paper 
\cite{ga} by Martin Gardner was the starting point 
of the popularity of the Bulgarian solitaire among mathematicians 
all over the world (see \cite{dr}).

\end{document}